\documentclass[12pt]{amsart}
\usepackage{amsmath}
\usepackage{amssymb}
\usepackage{bm}
\usepackage{graphicx}
\newtheorem{Theorem}{Theorem}
\newtheorem{Lemma}{Lemma}

\newtheorem{Corollary}{Corollary}

\begin{document}

\title[Non-localization for Graph Eigenfunctions]{Non-localization of Eigenfunctions on Large Regular Graphs}
\author{Shimon Brooks and Elon Lindenstrauss}
\thanks{E.L. was
supported in part by NSF grants DMS-0554345 and DMS-0800345.}

\maketitle

{\em Abstract: } We give a delocalization estimate for eigenfunctions of the discrete Laplacian on large $d+1$-regular graphs, showing that any subset of the graph supporting $\epsilon$ of the $L^2$ mass of an eigenfunction must be large.  For graphs satisfying a mild girth-like condition, this bound will be exponential in the size of the graph.

\section{Introduction}

The extent to which eigenfunctions of the Laplacian can localize in small sets is a question that has attracted much recent attention.  In the case of compact negatively curved surfaces, it is conjectured (see, eg.,  \cite{IwaniecSarnak}) that the sup-norms of Laplace eigenfunctions satisfy strong bounds (subexponential in the eigenvalue).  Another notion of delocalization is the Quantum Unique Ergodicity property, which roughly states that eigenfunctions become equidistributed in the ambient space. This property has been conjectured by Rudnick and Sarnak \cite{RS} to hold for all compact manifolds of negative sectional curvature, and  was proved for certain arithmetic manifolds (and a natural but specific choice of eigenbasis for the eigenfunctions) in  \cite{Lin}, \cite{Lior-Akshay} and \cite{Sound}.  Other work (eg., \cite{AN}) proves entropy estimates for eigenfunctions, a weaker notion of the inability to localize in small sets.

In this paper, we investigate eigenfunctions of the discrete Laplacian on large $d+1$-regular graphs, satisfying a mild condition (\ref{condition}) below (essentially asking that there not be too many short cycles through the same point).  We will use methods akin to the quantum chaos tools of \cite{AN} and others, and our results hold for {all} graphs satisfying (\ref{condition}) and any eigenfunction. It is likely much stronger results can be proven for random graphs. An earlier version of our result was also written in \cite{shimonthesis}.

To state the condition on our graphs, we let $\mathcal{T}_{d+1}$ be the $d+1$-regular tree (the universal cover of our $d+1$-regular graphs). For $f\in L^2(\mathcal{T}_{d+1})$, set  
$$\tilde{S}_{n}(f)(x) = d^{-n/2}\sum_{d(x,y)=n} f(y)$$
to be a normalized average over the sphere of radius $n$, and write $S_n$ for the projection of this operator to $L^2(\mathcal{G})$.  We assume that there exist $C$ and $\alpha>0$, such that the matrix coefficients of $S_n$ satisfy
\begin{equation}\label{condition}
\sup_{x\in\mathcal{G}}|| S_n\delta_x||_\infty \leq Cd^{-\alpha n} \text{ for all } n\leq N
\end{equation}
For example, if we set $N$ to be the radius of injectivity of the graph, then we may take $C=1$ and $\alpha = 1/2$.  Ideally, we would like to have $N\gtrsim \log{|\mathcal{G}|}$ (see below); at the very least, we should ensure that $N$ is large relative to the other parameters, tending to $\infty$ with $|\mathcal{G}|$.  

For random graphs, it is shown in \cite{MWW} that a large $d+1$-regular graph $\mathcal{G}$ almost surely does not have $2$ cycles of length at most $\left(\frac{1}{4}-\epsilon\right)\log_d{|\mathcal{G}|}$ that share an edge; in particular, this means that the condition (\ref{condition}) holds for almost all graphs.

More generally, we will assume control over the norm of $S_n$ as an operator from $L^p(\mathcal{G})$ to $L^q(\mathcal{G})$, for {\em some} conjugate pair $1\leq p <2 <q \leq \infty$ (see (\ref{hypothesis}) below).  This is equivalent to condition (\ref{condition}), up to changing the parameters $C$ and $\alpha$, but can give a better value for the bound $\delta$ in Theorem~\ref{main result}.

\begin{Theorem}\label{main result}
Let $\epsilon>0$, and $\mathcal{G}$ a $d+1$-regular graph satisfying 
\begin{equation}\label{hypothesis}
||S_n||_{L^p(\mathcal{G})\to L^q(\mathcal{G})} \leq C d^{-\alpha n} \text{ for all } n\leq N
\end{equation}
as an operator from $L^p(\mathcal{G})$ to $L^q(\mathcal{G})$, for some conjugate $1\leq p <2$ and $2<q\leq\infty$ (i.e., satisfying $\frac{1}{p}+\frac{1}{q}=1$).  
Then for any $L^2$-normalized eigenfunction $\phi$ on $\mathcal{G}$, any subset $E\subset \mathcal{G}$ satisfying 
$$\sum_{x\in E}|\phi (x)|^2 > \epsilon$$ 
must be of size 
$$|E| \gtrsim  d^{\delta N}$$
as $N\to\infty$, where $\delta = \delta(\epsilon, \alpha, p)$ can be taken to be $\delta = 2^{-7}\frac{\alpha p}{(2-p)}\epsilon^2$.
The implied constant depends on all parameters except $N$; namely $d$, $C$, $\alpha$, $p$, and $\epsilon$.
\end{Theorem}
Note that if $N\gtrsim \log_d{|\mathcal{G}|}$, then the conclusion of Theorem~\ref{main result} states that $|E|\gtrsim |\mathcal{G}|^{\delta'}$.

\section{Some Harmonic Analysis on the $d+1$-Regular Tree}
Throughout, we set $\mathcal{G}$ to be a $d+1$-regular graph.  We have the symmetric operator
$$T_df(x) = \frac{1}{\sqrt{d}}\sum_{d(x,y)=1}f(y)$$
and an orthonormal basis $\{\phi_j\}_{j=1}^{|\mathcal{G}|}$ of $L^2(\mathcal{G})$ consisting of $T_d$-eigenfunctions\footnote{The operator $T_d$ matches the operator defined above as $S_1$, though here we wish to emphasize the degree $d$ rather than the radius of the sphere.}.  (The discrete Laplacian on $\mathcal{G}$ can be written as $\Delta f = \left(\frac{\sqrt{d}}{d+1}T_d  - 1\right)f$, and so the eigenfunctions of $T_d$ are exactly the eigenfunctions of $\Delta$.)  

The universal cover of $\mathcal{G}$ is the $d+1$-regular tree, denoted $\mathcal{T}_{d+1}$.
Harmonic analysis on $\mathcal{T}_{d+1}$ has been well studied, see eg. \cite{FTP}.  For every $\lambda\in [-\frac{d+1}{\sqrt{d}},\frac{d+1}{\sqrt{d}}]$, there exists a unique {\bf spherical function} $\phi_\lambda$ satisfying:
\begin{itemize}
\item{$T_d\phi_\lambda = \lambda\phi_\lambda$.}
\item{$\phi_\lambda$ is radial; i.e., $\phi_\lambda(x) = \phi_\lambda(|x|)$ for all $x\in\mathcal{T}_{d+1}$, where $|x|$ denotes the distance from $x$ to the origin in $\mathcal{T}_{d+1}$.}
\item{$\phi_\lambda(0) = 1$.}
\end{itemize}
The last condition is simply a convenient normalization.  

We distinguish two parts of this spectrum:  the {\bf tempered spectrum} is the interval $[-2,2]$, and the {\bf untempered spectrum} is the part lying outside this interval, i.e.  $\pm(2,\frac{d+1}{\sqrt{d}}]$.  We will find it convenient to parametrize the spectrum by $\lambda = 2\cos{\theta_\lambda}$, where:
\begin{itemize}
\item{$\theta_\lambda\in [0,\pi]$ for $\lambda$ tempered.}
\item{$i\theta_\lambda = r_\lambda\in (0,\log{\sqrt{d}})$ for $\lambda$ untempered and positive.}
\item{$i\theta_\lambda + i\pi = r_{-\lambda}$ for $\lambda$ untempered and negative.}
\end{itemize}
 In this parametrization, we can write the spherical functions explicitly as \cite{dad}:
$$\phi_\lambda(x) = d^{-|x|/2}\left(\frac{2}{d+1}\cos{|x|\theta_\lambda} + \frac{d-1}{d+1}\frac{\sin{(|x|+1)\theta_\lambda}}{\sin{\theta_\lambda}}\right)$$
It will be convenient to use the Chebyshev polynomials 
\begin{eqnarray*}
P_n(\cos{\theta}) &  = & \cos{n\theta}\\
Q_n(\cos{\theta}) & = & \frac{\sin{(n+1)\theta}}{\sin{\theta}}
\end{eqnarray*}
of the first and second kinds, respectively.  With this notation the spherical functions become
\begin{equation}\label{spherical}
\phi_\lambda(x) = d^{-|x|/2}\left(\frac{2}{d+1}P_{|x|}(\lambda/2) + \frac{d-1}{d+1}Q_{|x|}(\lambda/2)\right)
\end{equation}

For any compactly supported radial function $k=k(|x|)$ on $\mathcal{T}_{d+1}$, the {\bf spherical transform} of $k$, denoted $h_k$, is given by 
$$h_k(\lambda) = \sum_{x\in\mathcal{T}_{d+1}}k(x)\phi_\lambda(x) = k(0) + (d+1)\sum_{n=1}^\infty d^{n-1}k(n)\phi_\lambda(n)$$
for all $\lambda\in [-\frac{d+1}{\sqrt{d}}, \frac{d+1}{\sqrt{d}}]$ (the sum is actually finite since $k$ is compactly supported).  

We have the {\bf Plancherel measure} $dm$ on $[-2,2]$ inverting the spherical transform on the tempered spectrum, i.e. 
$$\int_0^\pi \phi_{\lambda}(x) dm(\theta_\lambda) = \delta_0(x)$$
where $\delta_0$ is the $\delta$ function at $0$ on $\mathcal{T}_{d+1}$ given by 
$$\delta_0(x) = \left\{ \begin{array}{ccc} 1 & \quad & x=0\\ 0 & \quad & x\neq 0 \end{array}\right.$$

The Plancherel measure is absolutely continuous (with respect to Lebesgue measure on the semi-circle) and symmetric about $\pi/2$, and so its Fourier series is of the form
$$\frac{dm}{d\theta} = \sum_{j=0}^\infty c_j \cos{2j\theta}$$
The Plancherel measure is then given explicitly by \cite[Theorem 4.1]{FTP}
$$\int_0^\pi \cos{(2n\theta)} dm = \frac{1-d}{2d^{n}}$$
for $n>0$ (it is clear directly from the definitions that $\int_0^\pi dm = 1$).  

The spectrum of $T_d$ on $L^2(\mathcal{G})$ is contained in $[-\frac{d+1}{\sqrt{d}}, \frac{d+1}{\sqrt{d}}]$, and again we distinguish between the tempered eigenvalues in $[-2,2]$ and the untempered eigenvalues outside this interval.  Eigenfunctions of $T_d$ are also eigenfunctions of convolution with radial kernels; in fact, for a ``point-pair invariant" $k(x,y) = k(d(x,y))$ on $\mathcal{G}\times \mathcal{G}$, the eigenvalue for $\phi_j$ under convolution with $k$ depends only on $\lambda_j$, and is given by the spherical transform $h_k(\lambda_j)$ \cite{TerrWall}.

\section{The Main Estimate}

Our result centers on the following estimate for matrix coefficients of $P_n(\frac{1}{2}T_d)$; recall that $P_n(\cos{\theta}) = \cos{n\theta}$ are the Chebyshev polynomials (of the first kind).
\begin{Lemma}\label{estimate}
Let $\delta_0$ be the $\delta$-function supported at $0\in\mathcal{T}_{d+1}$, and $n$ a positive even integer.  Then
\begin{eqnarray*}
P_n(T_d/2)\delta_0(x) & = & \left\{ \begin{array}{ccc} 0  & \quad & |x|  \text{ odd } \quad \text{or} \quad |x|>n\\ \frac{1-d}{2d^{n/2}} & \quad & |x|<n \quad \text{and} \quad |x| \text{ even } \\ \frac{1}{2d^{n/2}} & \quad & |x| = n \end{array}\right.
\end{eqnarray*}
In particular, we have
$$P_n(T_d/2)\delta_0(x) \lesssim d^{-n/2}$$
\end{Lemma}

{\em Proof:}  Write $\delta_0 = \int_0^\pi \phi_{\lambda}dm(\theta_\lambda)$.  Then
since $\frac{1}{2}T_d\phi_\lambda = \cos{\theta_\lambda}\phi_\lambda$, we have
\begin{eqnarray*}
P_n(T_d/2)\delta_0(x) & = & \int_0^\pi (\cos{n\theta_\lambda)}\phi_\lambda(x) dm(\theta_\lambda)\\
& = & d^{-|x|/2}\int_0^{\pi}(\cos{n\theta_\lambda})
\left(\frac{2}{d+1}P_{|x|}(\lambda/2) + \frac{d-1}{d+1}Q_{|x|}(\lambda/2)\right)dm(\theta_\lambda)
\end{eqnarray*}
by substituting (\ref{spherical}) for the spherical functions.

Now since $n$ is even, both $\cos{n\theta_\lambda}$ and the Plancherel measure are  symmetric about $\pi/2$.  But if $|x|$ is odd, then both
\begin{eqnarray*}
P_{|x|}(\lambda/2) & = & \cos{|x|\theta_\lambda}\\
Q_{|x|}(\lambda/2) & = & \frac{\sin{(|x|+1)\theta_\lambda}}{\sin{\theta_\lambda}} = \cos{|x|\theta_\lambda} + \cos{\theta_\lambda}\frac{\sin{|x|\theta_\lambda}}{\sin{\theta_\lambda}}\\
& = & \cos{|x|\theta_\lambda} 
+ \cos{\theta_\lambda}\left( 1 + 2\sum_{j=1}^{(|x|-1)/2}\cos{2j\theta_\lambda}\right)
\end{eqnarray*}
are odd functions with respect to $\pi/2$.  Therefore the integral from $0$ to $\pi/2$ cancels with the integral from $\pi/2$ to $\pi$, and $P_n(T_d/2)\delta_0(x)$ vanishes for $|x|$ odd.

Now consider $|x|$ even, in which case we can write
$$Q_{|x|}(\lambda/2) = \frac{\sin{(|x|+1)\theta_\lambda}}{\sin{\theta_\lambda}} 
= 1 + 2\sum_{j=1}^{|x|/2}\cos{2j\theta_\lambda}$$
We will also make repeated use of the identity 
$$ 2\cos{\alpha}\cos{\beta} = \cos{(\alpha+\beta)} + \cos{(\alpha - \beta)}$$
If $|x|>n$, then we have
\begin{eqnarray*}
\lefteqn{d^{|x|/2}P_n(T_d/2)\delta_0(x)}\\ & = & \int_0^{\pi} (\cos{n\theta_\lambda}) 
\left(\frac{2}{d+1}P_{|x|}(\lambda/2) 
+ \frac{d-1}{d+1}Q_{|x|}(\lambda/2)\right)dm(\theta_\lambda)
\end{eqnarray*}
The left part of the integral yields
\begin{eqnarray*}
\lefteqn{\frac{2}{d+1}\int_0^{\pi} \cos{n\theta_\lambda}\cos{|x|\theta_\lambda}dm(\theta_\lambda)}\\ 
& = &
\frac{1}{d+1}\left(\int_0^{\pi}\cos{(|x|-n)\theta_\lambda}dm(\theta_\lambda) + \int_0^\pi \cos{(|x|+n)\theta_\lambda}dm(\theta_\lambda)\right)\\
& = & \frac{1}{d+1} \left(\frac{1-d}{2d^{(|x|-n)/2}} + \frac{1-d}{2d^{(|x|+n)/2}}\right)\\
& = & \frac{1-d}{d+1}\left(\frac{1}{2d^{(|x|-n)/2}}+\frac{1}{2d^{(|x|+n)/2}}\right)
\end{eqnarray*}
The right part, on the other hand, is
\begin{eqnarray*}
\lefteqn{\frac{d-1}{d+1}\int_0^{\pi}(\cos{n\theta_\lambda})
\left(1+2\sum_{j=1}^{|x|/2}\cos{2j\theta_\lambda}\right)dm(\theta_\lambda)}\\
& = & \frac{d-1}{d+1}\int_0^\pi \left(\cos{n\theta_\lambda} + 
\sum_{j=1}^{|x|/2}[\cos{(n+2j)\theta_\lambda} + \cos{(n-2j)\theta_\lambda}]\right)dm(\theta_\lambda)\\
& = & \frac{d-1}{d+1}\sum_{j=-|x|/2}^{|x|/2}\int_0^\pi\cos{(n+2j)\theta_\lambda}dm(\theta_\lambda)\\
& = & \frac{d-1}{d+1}\left(\sum_{j=-|x|/2}^{n/2-1}\frac{1-d}{2d^{n/2-j}}
+ 1 + \sum_{j=n/2+1}^{|x|/2}\frac{1-d}{2d^{j-n/2}}\right)\\
& = & \frac{d-1}{d+1}\left(\frac{1}{2d^{(n+|x|)/2}}
+ \frac{1}{2d^{(|x|-n)/2}}\right)
\end{eqnarray*}
since the sum telescopes.  Putting the two halves together gives
$P_n(T_d/2)\delta_0(x) = 0$ for all $|x|>n$.

If, however, $|x|\leq n$, then the sum does not quite telescope as before.  For $|x|<n$ we have instead
\begin{eqnarray*}
\lefteqn{d^{|x|/2}P_n(T_d/2)\delta_0(x)}\\
& = & \frac{1}{d+1}\left(\frac{1-d}{2d^{(n-|x|)/2}} + \frac{1-d}{2d^{(n+|x|)/2}}\right)
+ \frac{d-1}{d+1}\sum_{j=(n-|x|)/2}^{(n+|x|)/2}\frac{1-d}{2d^j}\\
& = & \frac{1-d}{d+1}\left(\frac{1}{2d^{(n-|x|)/2}} + \frac{1}{2d^{(n+|x|)/2}}\right)
+ \frac{d-1}{d+1}\left(-\frac{d}{2d^{(n-|x|)/2}}+\frac{1}{2d^{(n+|x|)/2}}\right)\\
& = & \frac{1-d}{d+1}\left(\frac{1}{2d^{(n-|x|)/2}}+\frac{d}{2d^{(n-|x|)/2}}\right)\\
& = & \frac{1-d}{2d^{(n-|x|)/2}}
\end{eqnarray*}
which finally gives 
$$P_n(T_d/2)\delta_0(x) = d^{-|x|/2}\frac{1-d}{2d^{(n-|x|)/2}} = \frac{1-d}{2d^{n/2}}$$
for $|x|<n$ and even.

If $|x|=n$, then we replace $\frac{1-d}{2d^{(n-|x|)/2}}$ with $1$, and the above calculation becomes
\begin{eqnarray*}
\lefteqn{d^{|x|/2}P_n(T_d/2)\delta_0(x)}\\
& = & \frac{1}{d+1}\left(1+ \frac{1-d}{2d^{(n+|x|)/2}} \right)
+ \frac{d-1}{d+1}\left( 1+ \sum_{j=1}^{(n+|x|)/2}\frac{1-d}{2d^j}\right)\\
& = & \frac{1}{d+1}\left(1 + \frac{1-d}{2d^{(n+|x|)/2}} \right)
+ \frac{d-1}{d+1}\left(\frac{1}{2} + \frac{1}{2d^{(n+|x|)/2}}\right)\\
& = & \frac{1}{d+1} + \frac{d-1}{2(d+1)}\\
& = & \frac{d+1}{2(d+1)} = 1/2
\end{eqnarray*}
and so $P_n(T_d/2)\delta_0(x) = \frac{1}{2d^{|x|/2}} = \frac{1}{2d^{n/2}}$ for $|x|=n$, as required.  $\Box$

\begin{Corollary}\label{p to q}
If $\mathcal{G}$ satisfies the hypothesis (\ref{hypothesis}) of Theorem~\ref{main result}, then
$$\left|\left| P_n\left(\frac{1}{2}T_d\right)\right|\right|_{L^p(\mathcal{G})\to L^q(\mathcal{G})} \lesssim d^{-\alpha n} $$
for all even positive integers $n\leq N$.
\end{Corollary}

{\em Proof:}
Thanks to Lemma~\ref{estimate}, we have that
$$ P_n(T_p/2) = \sum_{j=0}^{n/2-1} \frac{1-d}{2d^{n/2}} d^{j} S_{2j} + \frac{1}{2}S_n$$
as operators on functions of $\mathcal{G}$.  Therefore,
\begin{eqnarray*}
\left|\left| P_n\left(\frac{1}{2}T_d\right)\right|\right|_{L^p\to L^q} & \leq & 
\sum_{j=0}^{n/2-1} \frac{d-1}{2d^{n/2}} ||d^jS_{2j}||_{L^p\to L^q} + \frac{1}{2} ||S_n||_{L^p\to L^q}\\
& \leq & d^{-n/2+1} \sum_{j=0}^{n/2} d^j||S_{2j}||\\
& \leq & d^{-n/2+1} \sum_{j=0}^{n/2} Cd^{j(1-2\alpha)}\\
& \lesssim_{d,\alpha} &  d^{-n/2} \cdot C d^{n/2 - \alpha n} \lesssim d^{-\alpha n}
\end{eqnarray*}
as required.
 $\Box$

\section{Estimating the Mass of Small Sets}

Now we turn to the proof of Theorem~\ref{main result}.  We first consider tempered eigenfunctions; it will become clear how the argument is applied to untempered eigenfunctions as well.

\begin{Lemma}\label{Dirichlet}
Let $\epsilon>0$.  For any $\theta_0\in [0,\pi]$, there exists a kernel $k_{0}$ on $\mathcal{T}_{d+1}$ such that:
\begin{itemize}
\item{$k_0$ is supported on a ball of radius $N$.}
\item{The operator of convolution with $k_{0}$, call it $K_{\theta_0}(f) = f\ast k_{0}$, is bounded by 
$$||K_{\theta_0}||_{L^p(\mathcal{G})\to L^q(\mathcal{G})} \lesssim_{d,\alpha} C d^{-\frac{1}{128}\alpha N\epsilon^2}$$
as an operator from $L^p(\mathcal{G})$ to $L^q(\mathcal{G})$, where $C$ and $\alpha$ are the parameters of the hypothesis (\ref{hypothesis}).}
\item{The spherical transform of $k_{0}$ satisfies $h_{k_0 }\geq -1$ everywhere, and $h_{k_0}(\theta_0) > \epsilon^{-1}$.}
\end{itemize}
\end{Lemma}

{\em Proof:}  Set $M=\lfloor \epsilon^{-1}\rfloor$ and $R=\lceil \frac{1}{8}N\epsilon \rceil$ (one should think of $N$ as being much larger than $\epsilon^{-1}$, so that $R$ is large).  By Dirichlet's Theorem, we can find a positive integer $r\leq R$ such that $|r\theta_0 \mod{2\pi}| < 2\pi R^{-1} \leq \frac{2\pi}{N\epsilon}$.  There exists an even multiple of $r$, say $r' = 2lr$, such that $\frac{1}{16}R\epsilon \leq r' \leq 2R$ (if $r\geq \frac{1}{32}R\epsilon$, we can simply take $l=1$; otherwise there is a multiple of $r$ between $\frac{1}{32}R\epsilon$ and $\frac{1}{16}R\epsilon$, so take twice that multiple).  Moreover, since we can choose $2l\leq \frac{1}{16}R\epsilon$, we have $|r'\theta_0 \mod{2\pi}| < \frac{1}{8}\pi\epsilon \leq \frac{\pi}{8M}$.

We now set the spherical transform of $k_0$ to be $h_{k_0}(\theta) = F_{2M}(r'\theta) - 1$, where $F_{2M}$ is the Fej\`er kernel of order $2M$.  Since $r'\theta_0 \mod{2\pi}\in \left[-\frac{\pi}{8M}, \frac{\pi}{8M}\right]$ is close enough to $0$, we have 
$$F_{2M}(r'\theta_0) = \frac{1}{2M} \frac{\sin^{2}(2Mr'\theta_0)}{\sin^2(r'\theta_0)} > M+2$$ 
as long as $M\geq 4$, and therefore the eigenvalue of $\phi_{2\cos{\theta_0}}$ under convolution with $k_0$ will be $>M+1\geq \epsilon^{-1}$.  Moreover, since $F_{2M}$ is positive, the spherical transform of $k_0$ is bounded below by $-1$.  It remains to check the first two properties.

Now, by Corollary~\ref{p to q} of the main estimate, we see that the kernel whose spherical transform is $\cos{2j\theta}$--- i.e., the kernel of $P_{2j}(\frac{1}{2}T_d)$--- has norm $\lesssim d^{-2\alpha j}$ as a convolution operator from $L^p(\mathcal{G})$ to $L^q(\mathcal{G})$.  The spherical transform of $k_0$ is a sum of terms of the form $\frac{2M-j}{M}\cos{jr'\theta}$, where $j=1,2 ,\ldots, 2M$ (note that we eliminated the $j=0$ term by subtracting off the constant contribution to $F_{2M}$) and $r'\in 2\mathbb{Z}$.  Thus 
$$ ||K_{\theta_0}||_{L^p(\mathcal{G})\to L^q(\mathcal{G})} \lesssim  \sum_{j=1}^M d^{-\alpha jr'}
\lesssim_{d, \alpha} d^{-\alpha r'}$$
Then, since 
$$r'\geq \frac{1}{16}R\epsilon \geq \frac{1}{128}N\epsilon^2$$
this concludes the proof of Lemma~\ref{Dirichlet}.  $\Box$

We now wish to apply this convolution operator to examine the localization of eigenfunctions in small sets.  

{\em Proof of Theorem~\ref{main result}:}  Pick an eigenfunction $\phi_j$ of eigenvalue $\lambda_j$, and a set $E$ satisfying 
\begin{equation}\label{E}
||\phi_j||_{L^2(E)}^2 = ||\phi_j1_E||_{L^2(\mathcal{G})}^2 \geq \epsilon
\end{equation}
Define the operator $K_{j} = K_{\theta_{\lambda_j}}$ of Lemma~\ref{Dirichlet} (corresponding to $\theta_0=\theta_{\lambda_j}$) if $\lambda_j$ is tempered, or $K_j = K_{\theta_0=0}$ if $\lambda_j$ is untempered.  Observe that in either case $K_j$ satisfies
\begin{eqnarray}
\big|\langle K_{j}(\phi_j1_E),\phi_j1_E \rangle\big| & \leq & ||K_{j}(\phi_j1_E)||_q ||\phi_j1_E||_p \nonumber\\
& \leq & ||K_{j}||_{L^p\to L^q}||\phi_j1_E||_p^2\nonumber
\end{eqnarray}
By H\"older's Inequality, we have
$$ ||\phi_j1_E||_p^2 = ||\phi_j^p 1_E||_1^{2/p} \leq ||\phi_j^p||_{2/p}^{2/p}||1_E||_{2/(2-p)}^{2/p} = ||\phi_j||_2^2 \cdot |E|^{\frac{2-p}{p}}$$
so that
\begin{eqnarray}
\big|\langle K_{j}(\phi_j1_E),\phi_j1_E \rangle\big| & \leq & ||K_{j}||_{L^p\to L^q}||\phi_j1_E||_p^2\nonumber\\
& \leq & ||K_{j}||_{L^p\to L^q}||\phi_j||_2^2 \cdot|E|^{\frac{2-p}{p}}\nonumber\\
& \leq & ||K_{j}||_{L^p\to L^q} \cdot|E|^{\frac{2-p}{p}}\nonumber\\
& \lesssim_{d, \alpha} & C d^{-2^{-7}\alpha N \epsilon^2}\cdot |E|^{\frac{2-p}{p}}\label{|E|}
\end{eqnarray}
by Lemma~\ref{Dirichlet}.

On the other hand, decompose $\phi_j1_E$ spectrally as
$$\phi_j1_E = \langle \phi_j1_E, \phi_j\rangle \phi_j + g_\text{temp} + g_\text{untemp}$$
where $g_\text{temp}$ and $g_\text{untemp}$ are the tempered and untempered components of $\phi_j1_E$, respectively, excluding the $\phi_j$ component.
Notice that since 
$$|\langle \phi_j1_E, \phi_j\rangle|  =   ||\phi_j1_E||_2^2 \geq \epsilon$$
we have 
\begin{eqnarray}
||g_\text{temp}||_2^2 & \leq & ||\phi_j1_E||_2^2 - |\langle \phi_j1_E, \phi_j\rangle|^2
\nonumber\\
& = & ||\phi_j1_E||_2^2(1-||\phi_j1_E||_2^2)\nonumber\\
& \leq & ||\phi_j1_E||_2^2 (1-\epsilon) \label{gtemp}
\end{eqnarray}
Now the $K_{j}$-eigenvalue of any tempered eigenfunction is at least $-1$ by Lemma~\ref{Dirichlet}, and $K_{j}$ must be positive on the untempered eigenfunctions, since each term in the Fourier expansion of the Fejer kernel is of the form $\cos{2j\theta} = \cos{-i2jr}=\cosh({-2jr})>0$ for $\theta=-ir$, and similarly $\cos(2j\theta) = \cos(-i2jr-2j\pi) = \cosh({-2jr})$ for $\theta=-ir-\pi$.  Therefore
\begin{eqnarray}
\langle K_{j}(\phi_j1_E),\phi_j1_E\rangle & \geq & |\langle \phi_j1_E,\phi_j\rangle|^2 \langle K_{j}\phi_j,\phi_j\rangle - ||g_\text{temp}||_2^2\nonumber\\
& \geq & ||\phi_j1_E||_2^2\bigg(||\phi_j1_E||_2^2\langle K_{j}\phi_j,\phi_j\rangle - (1-\epsilon)\bigg)\label{spectral}
\end{eqnarray}
If $\lambda_j$ is tempered, then Lemma~\ref{Dirichlet} implies that the $\phi_j$-eigenvalue of $K_j$ is at least $\epsilon^{-1}$, whereby 
$$\langle K_j\phi_j,\phi_j\rangle \geq \epsilon^{-1}||\phi_j||_2^2 = \epsilon^{-1}$$
If $\lambda_j$ is untempered, then because $\cosh(2j\theta_j) > 1 = \cos(2j(0))$ we chose to use the same kernel as $\theta_0=0$ from the tempered case, and get that 
$$\langle K_{\theta_0=0}\phi_j,\phi_j\rangle \geq (F_{2M}(0)-1)||\phi_j||_2^2 = 2M-1 > \epsilon^{-1}$$
in the untempered case as well.  Applying (\ref{E}), we get from (\ref{spectral}) that
\begin{eqnarray}
\langle K_{j}(\phi_j1_E),\phi_j1_E\rangle
& \geq & ||\phi_j1_E||_2^2(||\phi_j1_E||_2^2\cdot \epsilon^{-1} - 1 +\epsilon)\nonumber\\
& \geq  & ||\phi_j1_E||_2^2(1 -1 + \epsilon) \nonumber\\
& \geq & \epsilon(\epsilon) = \epsilon^2\label{spectral side}
\end{eqnarray}
and combining (\ref{|E|}) with (\ref{spectral side}) yields
$$|E|^{\frac{2-p}{p}} \gtrsim_{d,\alpha} C^{-1}\epsilon^2 d^{2^{-7}\alpha\epsilon^2 N} \gtrsim d^{2^{-7}\alpha\epsilon^2 N}$$
which gives the bound of Theorem~\ref{main result}. $\Box$

\end{document}